\documentclass[12pt,a4paper]{article}

\usepackage[marked,reporte
]{myart98}
\ppnum{funct-an/9712003}{1997}{}
\def\ifundefined#1{\expandafter\ifx\csname#1\endcsname\relax}
\newcommand{\comment}[1]{}

\usepackage{amsfonts}

\newcommand{\algebra}[1]{\ensuremath{\mathfrak{#1}}}

\newcommand{\Cliff}[2][\comment]{\ensuremath{{\cal C}\kern-0.18em\ell(#1,#2)}}
\newcommand{\object}[2][\,]{\ensuremath{\mathrm{#2}#1}}
\newcommand{\Space}[2]{\ensuremath{ {\mathbb{#1}^{#2}} }}
\newcommand{\such}{\,\mid\,}

\newcommand{\FSpace}[2]{{\ensuremath{ #1_{#2} }}}

\ifundefined{qed}
    \DeclareMathSymbol{\qed}{0}{AMSa}{"03}
\fi
\newcommand{\norm}[1]{\left\| #1 \right\|}
\newcommand{\modulus}[1]{\left| #1 \right|}
\newcommand{\scalar}[2]{\langle #1,#2\rangle}

\providecommand{\eqref}[1]{\textup{(\ref{#1})}}

\newcommand{\person}[1]{\textsc{#1}}

\ifundefined{eoe}
   \DeclareMathSymbol{\eoe}      {\mathord}{AMSa}{"06}
\fi
\newcommand{\mod}{\mathrm{mod}\,}
     \theorembodyfont{\upshape}
\newtheorem{examplea}{\addtocounter{thm}{1}Example}
\newtheorem{exampleb}{Example}

\newcommand{\SL}{SL(2,\Space{R}{})}
\newcommand{\oper}[1]{\mathcal{#1}}
\newcommand{\matr}[4]{{\ensuremath{ \left( \begin{array}{cc}
#1 & #2 \\ #3 & #4 \end{array}\right) }}}
\newcommand{\vecbf}[1]{\mathbf{#1}}
\newcommand{\n}[1]{\mathsf{#1}}
\newcommand{\TSpace}[2]{\ensuremath{ { \widetilde{\mathbb{#1}}^{#2}} }}

\titleshort{Analysis in $\Space{R}{1,1}$}
\authorshort{Vladimir V. Kisil}

\begin{document}
\title{Analysis in $\Space{R}{1,1}$ \\
or the Principal Function Theory\thanks{Supported
by grant 3GP03196 of the FWO-Vlaanderen (Fund of Scientific 
Research-Flanders),
Scientific Research Network ``Fundamental Methods and Technique in
Mathematics'' and grant INTAS 93--0322--Ext.}}
\author{Vladimir V. Kisil\\
\normalsize 	Institute of Mathematics,
\normalsize                Economics and Mechanics,\\
\normalsize                Odessa State University, 
\normalsize                ul. Petra Velikogo, 2,\\
\normalsize                Odessa-57, 270057, UKRAINE\\
\normalsize                E-mail: \texttt{kisilv@member.ams.org}}
\maketitle
\begin{abstract}
We explore a function theory connected with the principal series
representation of $\SL$ in contrast to standard complex 
analysis connected with the discrete series.
We construct counterparts for the Cauchy integral formula, the Hardy space, 
the Cauchy-Riemann equation and the Taylor expansion.
\AMSMSC{30G35}{22E46}
\end{abstract}
\newpage
{\small \tableofcontents }
\section{Introduction}
\epigraph{You should complete your own original research in order to 
learn when it was done before.}{}{}
Connections between complex analysis (one variable, several complex 
variables, Clifford analysis) and its symmetry groups are known from its 
earliest days and are an obligatory part of the textbook on the 
subject~\cite{Cnops94a}, \cite{DelSomSou92}, \cite[\S~1.4, 
\S~5.4]{GilbertMurray91}, \cite{Krantz82}, \cite[Chap.~2]{Rudin80}.
However ideas about fundamental role of symmetries in function
theories outlined in~\cite{GilbertKunzeTomas86,Koranyi72a} were not
incorporated in a working toolkit of researchers yet.

It was proposed in~\cite{Kisil96c} to distinguish essentially different
function theories by corresponding group of symmetries.  Such a
classification is needed because not all seemingly different function
theories are essentially different~\cite{Kisil95c}.  But it is also
important that the group approach gives a constructive way to develop
essentially different function theories~\cite{Kisil96c,Kisil97a}, as well
as outlines an alternative ground for functional calculi of
operators~\cite{Kisil95i}.  In the mentioned papers all given examples
consider only well-known function theories.  While rearranging of known
results is not completely useless there was an appeal to produce a new
function theory based on the described scheme.

The theorem proved in~\cite{Kisil96d} underlines the similarity between
structure of the group of M\"obius transformations in spaces $ \Space{R}{n}
$ and $ \Space{R}{pq} $. This generates a hope that there exists a
non empty function theory in $ \Space{R}{pq} $. We construct such a theory
in the present paper for the case of $ \Space{R}{1,1} $. Other new function
theories based on the same scheme will be described
elsewhere~\cite{CnopsKisil97a}.

The format of paper is as follows. In Section~\ref{se:preliminaries} we 
introduce basic notations and definitions. We construct two function 
theories---the standard complex analysis and a function theory in 
$\Space{R}{1,1}$---in Section~\ref{se:theories}. Our consideration is based
on two different series of representation of $\SL$: discrete and principal.
We deduce in their terms the Cauchy integral formula, the Hardy spaces, the
Cauchy-Riemann equation, the Taylor expansion and their counterparts for 
$\Space{R}{1,1}$. Finally we collect in Appendix~\ref{se:appendix} several 
facts, which we would like (however can not) to assume well known. It may 
be a good idea to look through the 
Appendixes~\ref{pt:lie-alg}--\ref{pt:principal} between the reading of
Sections~\ref{se:preliminaries} and~\ref{se:theories}. 
Finally Appendix~\ref{pt:o-problems} states few among many open problems.
Our examples will be rather lengthy thus their (not always obvious) ends 
will be indicated by the symbol $\eoe$.

There are no claims about novelty any particular formula for analysis in 
$\Space{R}{1,1}$. Most of them are probably known in some different form in
the theory of special functions (see for example~\ref{pr:taylor-b}). 
Moreover our general ideas are very close
to~\cite{GilbertKunzeTomas86}. However the composition of results as a 
function theory parallel to the complex analysis is believed to be new.

\section{Preliminaries} \label{se:preliminaries}
Let $\Space{R}{pq}$ be a real $n$-dimensional vector space, where
$n=p+q$ with a fixed frame $e_1$, $e_2$, \ldots, $e_p$, $e_{p+1}$,
\ldots, $e_n$ and with the nondegenerate bilinear
form $B(\cdot,\cdot)$ of the signature $(p,q)$, which is diagonal in the
frame $e_i$, i.e.:
\begin{displaymath}
B(e_i,e_j)=\epsilon_i \delta_{ij}, \textrm{ where }
\epsilon_i=\left\{\begin{array}{ll}
1, & i=1,\ldots,p\\
-1, & i=p+1,\ldots,n
\end{array}\right.
\end{displaymath}
and $\delta _{ij}$ is the Kronecker delta. In particular the usual
Euclidean space $\Space{R}{n}$ is $\Space{R}{0n}$. Let $\Cliff[p]{q}$ be
the \emph{real Clifford algebra} generated by $1$, $e_j$, $1\leq j\leq n$
and the relations 
\begin{displaymath}
e_i e_j + e_j e_i =-2B(e_i,e_j).
\end{displaymath}
We put $e_0=1$ also. Then there is the natural embedding
$\algebra{i}: \Space{R}{pq}\rightarrow
\Cliff[p]{q}$. We identify $\Space{R}{pq}$ with its image under
$\algebra{i}$ and
call its elements \emph{vectors}. There are two linear
anti{-}automorphisms $*$ (reversion) and $-$
(main anti{-}automorphisms) 
and automorphism $'$
of $\Cliff[p]{q}$ defined on its
basis $A_\nu=e_{j_1}e_{j_2}\cdots e_{j_r}$, $1\leq j_1 <\cdots<j_r\leq
n$ by the rule:
\begin{eqnarray*}
(A_\nu)^*= (-1)^{\frac{r(r-1)}{2}} A_\nu, \qquad
\bar{A}_\nu= (-1)^{\frac{r(r+1)}{2}} A_\nu,\qquad
A_\nu'= (-1)^{r} A_\nu.
\end{eqnarray*}
In particular, for vectors, $\bar{\vecbf{x}}=\vecbf{x}'=-\vecbf{x}$ and
$\vecbf{x}^*=\vecbf{x}$.

It is easy to see that $\vecbf{x}\vecbf{y}=\vecbf{y}\vecbf{x}=1$
for any $\vecbf{x}\in\Space{R}{pq}$ such that
$B(\vecbf{x},\vecbf{x})\neq 0$
and $\vecbf{y}={\bar{\vecbf{x}}}\,{\norm{\vecbf{x}}^{-2}}$, which is
the \emph{Kelvin inverse} of $\vecbf{x}$.
Finite products of invertible vectors are invertible in $\Cliff[p]{q}$
and form the \emph{Clifford group} $\Gamma(p,q)$. Elements
$a\in\Gamma(p,q)$ such that
$a\bar{a}=\pm 1$ form the $\object[(p,q)]{Pin}$ group---the double cover
of the group of orthogonal rotations $\object[(p,q)]{O}$. We also
consider~\cite[\S~5.2]{Cnops94a} $T(p,q)$ to be the set of all
products of vectors in $\Space{R}{pq} $.

Let $(a, b, c, d)$ be a quadruple from $T(p,q)$ with
the properties:
\begin{enumerate}
\item $(ad^*-bc^*)\in \Space{R}{}\setminus {0}$;
\item $a^*b$, $c^*d$, $ac^* $, $ bd^*$ are vectors.
\end{enumerate}
Then~\cite[Theorem~5.2.3]{Cnops94a}
$2\times 2$-matrixes \matr{a}{b}{c}{d} form the
group $\Gamma(p+1,q+1)$ under the
usual matrix multiplication. It has a representation
$\pi_{\Space{R}{pq} }$ \comment{we denote its restriction to any
subgroup by the same notation.} by transformations of
$\overline{\Space{R}{pq} }$ given by:
\begin{equation}\label{eq:sp-rep}
\pi_{\Space{R}{pq}}\matr{a}{b}{c}{d} :
\vecbf{x} \mapsto (a\vecbf{x}+b)(c\vecbf{x}+d)^{-1},
\end{equation}
which form the \emph{M\"obius} (or
the \emph{conformal}) group of $\overline{\Space{R}{pq}}$. 
Here $\overline{\Space{R}{pq}}$ the compactification of $\Space{R}{pq} $
by the ``necessary number 
of points'' (which form the light cone) at infinity 
(see~\cite[\S~5.1]{Cnops94a}).
The analogy with fractional-linear transformations of the complex line
\Space{C}{} is useful, as well as representations of shifts
$\vecbf{x}\mapsto \vecbf{x}+y$, orthogonal rotations
$\vecbf{x}\mapsto k(a)\vecbf{x}$, dilations
$\vecbf{x}\mapsto \lambda \vecbf{x}$, and the Kelvin inverse
$\vecbf{x}\mapsto \vecbf{x}^{-1}$ by the
matrixes \matr{1}{y}{0}{1}, \matr{a}{0}{0}{{a}^{*-1}},
\matr{\lambda^{1/2}}{0}{0}{\lambda^{-1/2}}, \matr{0}{-1}{1}{0}
respectively. We also use the agreement of~\cite{Cnops94a} that $ 
\frac{a}{b} $ always denotes $ab^{-1}$ for $a$, $b\in \Cliff[p]{q}$.

\section{Two Function Theories Associated with representations of $\SL$} 
\label{se:theories}
\subsection{Unit Disks in $ \Space{R}{0,2} $ and $ \Space{R}{1,1} $}
We describe a coherent states and wavelet transform connected with a 
homogeneous space $\Omega=G/H$ and unitary irreducible representation
$\pi$ of $G$, which is induced by a character of 
$H$~\cite[\S~13.2]{Kirillov76}. This representation is assumed to be a 
square integrable with respect to $ \Omega $ (see below).

Let $G$ be a locally compact group and $H$ be its closed subgroup.  Let $
\Omega=G / H$ and $s:  \Omega \rightarrow G$ be a continuous
mapping~\cite[\S~13.2]{Kirillov76}.  Then any $g\in G$ has a unique
decomposition of the form $g=s(\omega)h$, $\omega\in \Omega$ and we will
write $\omega=s^{-1}(g)$, $h=r(g)={(s^{-1}(g))}^{-1}g$.  Note that $ \Omega
$ is a left $G$-homogeneous space with an action defined in terms of $s$ as
follow:  $g:  a \mapsto g\cdot \omega=s^{-1}(g^{-1}* s(\omega)) $, where
$*$ is the multiplication on $G$.

The main example is provided by group 
$G=\SL$ (books~\cite{HoweTan92,Lang85,MTaylor86} are our
standard references about $ \SL $ and its
representations) consisting of $ 2\times 2$ matrices $
\matr{a}{b}{c}{d}$ with real entries and determinant $ad-bc=1$. 
See Appendix~\ref{pt:lie-alg} for description of its Lie algebra.
We will construct two series of examples. One is connected with discreet 
series representation and produces the core of standard complex analysis. 
The second will be its mirror in principal series representations and 
create parallel function theory. $ \SL $ has also other type 
representation, which can be of particular interest in other 
circumstances. However the discreet series and principal ones stay 
separately from others (in particular by being the support of the
Plancherel measure~\cite[\S~VIII.4]{Lang85},~\cite[Chap.~8,
(4.16)]{MTaylor86}) and are in a good resemblance each other.

\begin{examplea}
Via identities
\begin{displaymath}
\alpha= \frac{1}{2}(a+d-ic+ib), \qquad \beta= \frac{1}{2}(c+b -ia+id)
\end{displaymath}
we have isomorphism of $ \SL $ with group $SU(1,1)$ of $
2\times 2$ matrices with complex entries of the form $ \matr{ \alpha}{
\beta}{ \bar{\beta}}{\bar{ \alpha}}$ such that $ \modulus{ \alpha }^2 - 
\modulus{ \beta }^2=1 $. We will use the last form for  $ \SL$ for complex 
analysis in unit disk $\Space{D}{}$.

$ \SL $ has the only non-trivial compact closed subgroup
$K$, namely the group of matrices of the form
$h_{\psi}=\matr{e^{i\psi}}{0}{0}{e^{-i\psi}}$.
Now any $ g\in \SL $ has a unique decomposition of the
form 
\begin{eqnarray}
\matr{\alpha}{\beta}{\bar{\beta}}{\bar{\alpha}} 
& =& \modulus{\alpha} 
\matr{1}{\beta\bar{\alpha}^{-1}}{\bar{\beta}\alpha^{-1}}{1} \matr{ 
\frac{{\alpha}}{ \modulus{\alpha} } }{0}{0}{\frac{\bar{\alpha}}{ 
\modulus{\alpha} } }
\nonumber \\
& =& \frac{1}{\sqrt{1- \modulus{a}^2 }}
\matr{1}{a}{\bar{a}}{1} 
\matr{e^{i\psi}}{0}{0}{e^{-i\psi}} 
\end{eqnarray}
where $\psi=\Im \ln (\alpha)$, $a=\beta\bar{\alpha}^{-1}$, and $
\modulus{ a } < 1$ because $ \modulus{ \alpha }^2 - \modulus{ \beta
}^2=1 $. Thus we can identify $\SL / H$ with the unit
disk $ \Space{D}{} $ and define mapping $s: \Space{D}{} \rightarrow
\SL $ as follows
\begin{equation} \label{eq:def-s-a}
s: a \mapsto \frac{1}{\sqrt{1- \modulus{a}^2 }}
\matr{1}{a}{\bar{a}}{1}.
\end{equation}
Mapping $ r: G \rightarrow H$ associated to $s$  is
\begin{equation} \label{eq:def-r-a}
r: \matr{\alpha}{\beta}{\bar{\beta}}{\bar{\alpha}} \mapsto \matr{ 
\frac{{\alpha}}{
\modulus{\alpha} }}{0}{0}{\frac{\bar{\alpha}}{ \modulus{\alpha} }}
\end{equation}

The invariant measure $ d\mu(a) $ on $ \Space{D}{} $ coming from
decomposition $dg=d\mu(a)\, dk$, where $dg$ and $dk$ are Haar measures
on $G$ and $K$ respectively, is equal to 
\begin{equation} \label{eq:def-m-a}
d\mu(a)= \frac{da}{{(1- \modulus{a}^2)}^2}.
\end{equation}
with $da$---the standard Lebesgue measure on $ \Space{D}{} $.

The formula $g: a \mapsto g\cdot a=s^{-1}(g^{-1} * s(a)) $ associates with
a matrix $ g^{-1}=\matr{ \alpha}{ \beta}{ \bar{\beta}}{ \bar{\alpha}}$ the
fraction-linear transformation of $ \Space{D}{} $ of the form
\begin{equation} \label{eq:fr-lin-a}
g: z \mapsto g\cdot z=
\frac{\alpha z + \beta}{\bar{ \beta}z +\bar{\alpha}},  
\qquad g^{-1}=\matr{\alpha}{\beta}{\bar{\beta}}{\bar{\alpha}},
\end{equation}
which also can be considered as a transformation of $\dot{ \Space{C}{} }$
(the one-point compactification of $ \Space{C}{} $). $\eoe$
\end{examplea}

\begin{exampleb}
We will describe a version of previous formulas corresponding to 
geometry of unit disk in $ \Space{R}{1,1}$. For generators $e_1$ and $e_2$ 
of $ \Space{R}{1,1} $ (here $e_1^2=-e^2_2=-1$.) we see that matrices 
$\matr{a}{be_2}{ce_2}{d}$ again give a realization of $\SL$. Making 
composition with the Caley transform
\begin{displaymath}
T=\frac{1}{ {2} } \matr{1}{e_2}{e_2}{-1} \matr{1}{e_1}{e_1}{1} =
\frac{1}{ {2} } \matr{1+e_2e_1}{e_1+e_2}{e_2-e_1}{e_2e_1-1}
\end{displaymath} 
and its inverse 
\begin{displaymath}
T^{-1}=\frac{1}{ {2} } \matr{1}{-e_1}{-e_1}{1} \matr{1}{e_2}{e_2}{-1}=
\frac{1}{2} \matr{1-e_1e_2}{e_2+e_1}{e_2-e_1}{-1-e_1e_2}
\end{displaymath}
(see analogous calculation in~\cite[\S~IX.1]{Lang85}) we obtain another
realization of $\SL$:
\begin{equation} \label{eq:caley}
\frac{1}{4} \matr{1-e_1e_2}{e_2+e_1}{e_2-e_1}{-1-e_1e_2}
\matr{a}{be_2}{ce_2}{d} 
\matr{1+e_2e_1}{e_1+e_2}{e_2-e_1}{e_2e_1-1}
=\matr{\n{a}}{\n{b}}{\n{b}'}{\n{a}'},
\end{equation}                     
where                              
\begin{equation}                   
\n{a}= \frac{1}{2}(a(1-e_1e_2)+d(1+e_1e_2)), \qquad 
\n{b}= \frac{1}{2} (b(e_1-e_2)+c(e_1+e_2)).
\end{equation}
It is easy to check that the condition $ad-bc=1$ implies the following
value of the pseudodeterminant of the matrix $ \n{a}(\n{a}')^*
-\n{b}(\n{b}')^*= \n{a}\bar{\n{a}}-\n{b} \bar{\n{b}}=1$.  We also observe
that $\n{a}$ is an even Clifford number and $\n{b}$ is a vector thus
$\n{a}'=\n{a}$, $\n{b}'=-\n{b}$.

Now we consider the decomposition
\begin{equation}
\matr{\n{a}}{\n{b}}{-\n{b}}{\n{a}} = 
\modulus{\n{a}} \matr{1}{\n{b}\n{a}^{-1}}{-\n{b}\n{a}^{-1}}{1}   \matr{ 
\frac{\n{a}}{ \modulus{\n{a}}} }{0}{0}{\frac{\n{a}}{ \modulus{\n{a}}}}.
\end{equation}
It is seen directly, or alternatively follows from general characterization 
of $\Gamma(p+1,q+1)$~\cite[Theorem~5.2.3(b)]{Cnops94a}, that
$\n{b}\n{a}^{-1}\in \Space{R}{1,1}$.  Note that now we cannot derive from $
\n{a}\bar{\n{a}}-\n{b} \bar{\n{b}}=1$ that $\n{b}\n{a}^{-1}
\overline{\n{b}\n{a}^{-1}}= -{(\n{b}\n{a}^{-1})}^2 <1 $ because $
\n{a}\bar{\n{a}}$ can be positive or negative (but we are sure that
${(\n{b}\n{a}^{-1})}^2\neq -1$).  For this reason we cannot define the unit
disk in $ \Space{R}{1,1} $ by the condition $ \modulus{\n{u}}<1 $ in a way
consistent with its M\"obius transformations.  This topic will be discussed
elsewhere in more details and illustrations~\cite{CnopsKisil96a}.  We
describe a ready solution in Appendix~\ref{pt:disk}.

Matrices of the form 
\begin{displaymath}
\matr{\n{a}}{0}{0}{\n{a}'}=\matr{e^{e_1e_2\tau}}{0}{0}{e^{e_1e_2\tau}}, 
\qquad \n{a}={e^{e_1e_2\tau}}=\cosh\tau+e_1e_2\sinh\tau, \quad \tau\in 
\Space{R}{}
\end{displaymath}
comprise a subgroup of $\SL$ which we denote by $A$.  This subgroup is an
image of the subgroup $A$ in the Iwasawa decomposition
$\SL=ANK$~\cite[\S~III.1]{Lang85} under the
transformation~\eqref{eq:caley}.

We define an embedding $s$ of $\TSpace{D}{}$ for our realization of 
$\SL$ by the formula:
\begin{equation} \label{eq:def-s-b}
s: \n{u} \mapsto \frac{1}{ \sqrt[]{1+\n{u}^2} } \matr{1}{\n{u}}{-\n{u}}{1}.
\end{equation}
The formula $g: \n{u} \mapsto s^{-1}(g\cdot s(\n{u}))$ associated with a 
matrix $g^{-1}=\matr{\n{a}}{\n{b}}{-\n{b}}{\n{a}}$ gives the 
fraction-linear transformation $ \TSpace{D}{} \rightarrow \TSpace{D}{}$ of 
the form:
\begin{equation} \label{eq:fr-lin-b}
g: \n{u} \mapsto g\cdot \n{u}= \frac{\n{a}\n{u}+\n{b}}{-\n{b}\n{u}+\n{a}}, 
\qquad g^{-1}=\matr{\n{a}}{\n{b}}{-\n{b}}{\n{a}}
\end{equation}
The mapping $ r: G \rightarrow H$ associated to $s$ defined 
in~\eqref{eq:def-s-b} is
\begin{equation} \label{eq:def-r-b}
r: \matr{\n{a}}{\n{b}}{-\n{b}}{\n{a}} \mapsto \matr{ \frac{\n{a}}{ 
\modulus{\n{a}} }}{0}{0}{\frac{\n{a}}{ \modulus{\n{a}} }}
\end{equation}

And finally the invariant measure on $ \TSpace{D}{} $ 
\begin{equation} \label{eq:def-m-b}
d\mu(\n{u})= \frac{d\n{u}}{{(1+\n{u}^2)}^2} = \frac{du_1 
du_2}{{(1-u_1^2+u_2^2)}^2}.
\end{equation}
follows from the elegant consideration in~\cite[\S~6.1]{Cnops94a}. $\eoe$
\end{exampleb}
We hope the reader notes the explicit similarity between these two 
examples.
Following examples will explore it further.

\subsection{Reduced Wavelet Transform---the Cauchy Integral Formula}
Let $ \chi: H \rightarrow \Space{C}{} $ be a unitary character of $H$,
which induces in the sense of Mackey an irreducible unitary representation
$\pi$ of $G$~\cite[\S~13.2]{Kirillov76}. We denote by the same letter $\pi$
the canonical realization of this representation in the space $ 
\FSpace{L}{2}(\Omega) $ given by the
formula~\cite[\S~13.2.(7)--(9)]{Kirillov76}:
\begin{equation} \label{eq:def-ind}
[\pi(g) f](\omega)= \chi_0(r(g^{-1} * s(\omega)))  f(g\cdot \omega),
\qquad 
\chi_0(h)=\chi(h)\left( 
\frac{d\mu(h\cdot \omega)}{d\mu(\omega)} \right)^{ \frac{1}{2} },
\end{equation}
where $g\in G$, $\omega\in\Omega$, $h\in H$ and $r:  G \rightarrow H$, $s:
\Omega \rightarrow G$ are functions defined at the beginning of this
Section; $*$ denotes multiplication on $G$ and $\cdot$ denotes the action
of $G$ on $\Omega$ from the left.

Let $ \pi_0 $ be a representation of $G$ by operators on a 
Hilbert space $ \FSpace{L}{2}(X,d\nu) $ such that its restriction to a 
subspace $ \FSpace{F}{2}(X,d\nu) \subset \FSpace{L}{2}(X,d\nu)$ 
is an irreducible unitary representation unitary equivalent to $\pi$. 
It follows in particular that 
\begin{equation} \label{eq:h-act}
[\pi_0(h) f](x) =  \chi_0(h) f(h\cdot x), \qquad \forall h\in H, \quad f(x) 
\in
\FSpace{F}{2}(X,d\nu).
\end{equation}
\begin{defn}
Let $f_0(x) \in \FSpace{F}{2}(X,d\nu) $ be an eigenfunction for all 
$\pi_0(h)$, $h\in H$, i.e, $\pi_0(h) f_0(x)=\chi_0(h) f_0(x)$. It is called the
\emph{vacuum vector}. For a fixed vacuum vector $f_0(x)$ we define the
\emph{reduced wavelet
transform}~\cite{Kisil97a} $\oper{W}: \FSpace{L}{2}(X,d\nu) \rightarrow
\FSpace{L}{\infty}(\Omega)$ by the formula
\begin{equation} \label{eq:def-w}
[\oper{W} f] (\omega)= \scalar{f}{ \pi_0(s(\omega)) f_0}.
\end{equation}
We say that a representation $\pi_0$ is \emph{square integrable $\mod H$}
if the reduced wavelet transform $[\oper{W}f_0](\omega)$ of the vacuum
vector $f_0(x)$ is square integrable on $ \Omega $.
\end{defn}

The reduced wavelet transform has following properties.
\begin{prop}
The reduced wavelet transform is a continuous mapping $\oper{W}:
\FSpace{L}{2}(X,d\nu) \rightarrow \FSpace{L}{\infty}(\Omega)$.
The kernel of the reduced wavelet transform $\oper{W}$ is the orthogonal 
completion $ \FSpace{F}{2} (X)^\perp $ of $ \FSpace{F}{2}(X) $. 
If $\pi_0$ is square integrable $\mod H$ then 
\begin{enumerate}
\item $\oper{W}$ maps $ \FSpace{L}{2}(X) $ to $ \FSpace{L}{2}(\Omega) $.
\item The image of $\oper{W}$ is the closure of linear combinations of 
vectors $
\widehat{f}_\omega= \oper{W} \pi_0(s(\omega)) f_0 $.
\item  $\oper{W}$ intertwines $\pi_0$ and $\pi$: 
\begin{equation} \label{eq:w-inter}
\oper{W} \pi_0(g) = \pi(g) \oper{W}.
\end{equation}
\end{enumerate}
\end{prop}
\begin{proof}
It is obvious from~\eqref{eq:def-w} that every function in  $ \FSpace{F}{2} 
(X)^\perp $ belongs to the kernel of $\oper{W}$. The intersection of the 
kernel with  $ \FSpace{F}{2} (X)$ is a subspace of $\FSpace{F}{2} (X)$ 
invariant under $\pi_0$. Due to non-degeneracy of $\pi_0$ and its 
irreducibility on  $ \FSpace{F}{2} (X)$ this invariant subspace is trivial.

Due to the irreducibility of $\pi_0$ on  $ \FSpace{F}{2} (X)$ the vacuum 
vector $f_0(x)$ is cyclic in $ \FSpace{F}{2} (X)$.  Thus linear
combinations of $\pi_0(g) f_0(x)$, $g\in G$ are dense in $ \FSpace{F}{2}
(X)$.  Due to the trivial action~\eqref{eq:h-act} of $H$ the same is true
for a smaller set $\pi_0(s(\omega)) f_0$, $\omega\in\Omega$.  Thus
the $\oper{W} \pi_0(s(\omega)) f_0$, $\omega\in\Omega$ are dense in the 
image of $\oper{W}$.

Finally,~\eqref{eq:w-inter} follows from the direct calculation:
\begin{eqnarray*}
[\oper{W} \pi_0(g) f] (\omega) 
& = & \scalar{\pi_0(g)f}{\pi_0(s(\omega))f_0} \\
& = & \scalar{f}{\pi_0(g^{-1})\pi_0(s(\omega))f_0} \\
& = &  \scalar{f}{\pi_0(g^{-1} *s( \omega))f_0} \\
& = & \chi_0(r(g^{-1}*s(\omega))) \scalar{f}{\pi_0(s(g\cdot \omega))f_0} \\
& = & \chi_0(r(g^{-1}*s(\omega))) [\oper{W} f](g\cdot \omega) \\
&=& \pi(g) [\oper{W}f](\omega).
\end{eqnarray*}
\end{proof}
\begin{examplea}
We continue to consider the case of $G=\SL$ and $H=K$. The compact group
$K\sim \Space{T}{}$ has a discrete set of characters $ \chi_m(h_\phi)=
e^{-i m\phi} $, $m\in \Space{Z}{} $. We drop the trivial character $\chi_0$
and remark that characters $\chi_m$ and $ \chi_{-m} $ give similar
holomorphic and \emph{anti}holomorphic series of representations.  Thus we
will consider only characters $\chi_m$ with $m=1,2,3,\ldots$.

There is a difference in behavior of characters $\chi_1$ and $ \chi_m $ for 
$m=2,3,\ldots$ and we will consider them separately.

First we describe $\chi_1$. Let us take $X=\Space{T}{}$---the unit circle 
equipped with the standard Lebesgue measure $ d\phi $ normalized in such a
way that 
\begin{equation} \label{eq:n-measure}
\int_{ \Space{T}{} } \modulus{f_0(\phi)}^2\, d\phi=1 \textrm{ with }
f_0(\phi)\equiv 1. 
\end{equation}
From ~\eqref{eq:def-s-a} and~\eqref{eq:def-r-a} one can find that
\begin{displaymath}
r(g^{-1}*s(e^{i\phi}))= \frac{ \bar{\beta}e^{i\phi}+\bar{\alpha}}
{ \modulus{ \bar{\beta}e^{i\phi}+\bar{\alpha}} }, 
\qquad
g^{-1}=\matr{\alpha}{\beta}{ \bar{\beta} }{ \bar{\alpha} }.
\end{displaymath}
Then the action of $G$ on $\Space{T}{}$ defined by~\eqref{eq:fr-lin-a}, 
the equality  $ {d(g\cdot \phi)}/{d\phi}= \modulus{\bar{\beta} e^{i\phi} + 
\bar{\alpha} }^{-2} $  and the character $ \chi_1 $ give the following
realization of the formula~\eqref{eq:def-ind}: 
\begin{equation} \label{eq:g-transform}
[\pi_1(g) f](e^{i\phi})= \frac{1}{ \bar{\beta} e^{i\phi} + \bar{ \alpha }} 
f \left( \frac{  { \alpha }e^{i\phi}+{\beta}}{\bar{\beta} e^{i\phi} 
+ \bar{ \alpha }} \right).
\end{equation}
This is a unitary representation---the \emph{mock discrete 
series} of $ \SL$~\cite[\~8.4]{MTaylor86}. It is easily seen that $K$ acts 
in a trivial way~\eqref{eq:h-act} by multiplication by $\chi(e^{i\phi})$.
The function $ f_0(e^{i\phi})\equiv 1 $ mentioned in~\eqref{eq:n-measure}
transforms as follows
\begin{equation} \label{eq:t-vacuum}
[\pi_1(g) f_0](e^{i\phi})= \frac{1}{ \bar{\beta} e^{i\phi} + \bar{ \alpha 
}}
\end{equation}
and in particular has an obvious property
$[\pi_1(h_\psi)f_0](\phi)=e^{i\psi} f_0(\phi)$, i.e.  it is a \emph{vacuum
vector} with respect to the subgroup $H$.  The smallest linear subspace $
\FSpace{F}{2}(X) \in \FSpace{L}{2}(X) $ spanned by~\eqref{eq:t-vacuum}
consists of boundary values of analytic functions in the unit disk and is
the \emph{Hardy space}.  Now the reduced wavelet transform~\eqref{eq:def-w}
takes the form
\begin{eqnarray}
\widehat{f}(a)=[\oper{W} f] (a) & =  &
\scalar{f(x)}{\pi_1(s(a))f_0(x)}_{\FSpace{L}{2}(
X ) } \nonumber \\
& = & \int_{\Space{T}{}} f(e^{i\phi}) \frac{ \sqrt{ 1-\modulus{a}^2 }
}{ \overline{ 
\bar{a}e^{i\phi} - 1} }\,d\phi \nonumber \\
& = & \frac{\sqrt{ 1-\modulus{a}^2 }}{i}  \int_{\Space{T}{}} 
\frac{f(e^{i\phi})}{{a}-e^{i\phi}} ie^{i\phi}\,d\phi 
\nonumber \\
& = & \frac{\sqrt{ 1-\modulus{a}^2 }}{i}
 \int_{\Space{T}{}}  \frac{f(z)}{{a}-z}\,dz, \label{eq:cauchy}
\end{eqnarray}
where $z=e^{i\phi}$. Of course~\eqref{eq:cauchy} is the \emph{Cauchy
integral formula} up to factor ${2\pi } {\sqrt{ 1-\modulus{a}^2 }} $.
Thus we will write $f(a)={\left({2\pi \sqrt{ 1-\modulus{a}^2 
}}\right)}^{-1}
\widehat{f}(a) $ for analytic extension of $f(\phi)$ to the unit disk.
The factor $2\pi $ is due to our normalization~\eqref{eq:n-measure} and
$\sqrt{ 1-\modulus{a}^2 }$ is connected with the invariant measure on $
\Space{D}{} $. 

Let us now consider characters $\chi_m$ ($m=2,3,\ldots$). 
These characters together with action~\eqref{eq:fr-lin-a} of $G$ give following
realization of~\eqref{eq:def-ind}:
\begin{equation} \label{eq:rep-dis}
[\pi_m(g) f](w)=f\left(  \frac{\alpha w + \beta}{\bar{\beta} w + 
\bar{\alpha}} \right) {(\bar{\beta} w + \bar{\alpha} )}^{-m}.
\end{equation}
For any integer $m\geq 2$ one can select a measure 
\begin{displaymath}
d\mu_m(w)=4^{1-m}{(1- \modulus{w}^2 )}^{m-2} dw,
\end{displaymath}
where $dw$ is the standard Lebesgue measure on $ \Space{D}{} $, such
that~\eqref{eq:rep-dis} become unitary
representations~\cite[\S~IX.3]{Lang85},~\cite[\S~8.4]{MTaylor86}. 
These are discrete series.
 
If we again select $f_0(w)\equiv 1$ then
\begin{displaymath}
[\pi_m(g) f_0](w)= {(\bar{\beta} w + \bar{\alpha} )}^{-m}.
\end{displaymath}
In particular $[\pi_m(h_\phi) f_0](w)= e^{im\phi} f_0(w)$ so this  
again is a vacuum vector with respect to $K$. The irreducible subspace $ 
\FSpace{F}{2}( \Space{D}{} )$ generated by $f_0(w)$ consists of analytic
functions and  is the $m$-th Bergman space (actually \person{Bergman}
considered only $m=2$). Now the transformation~\eqref{eq:def-w} 
takes the form                                                                    
\begin{eqnarray*}
\widehat{f}(a) & = & \scalar{f(w)}{[\pi_m(s(a)) f_0](w)}\\
& = & {\left(\sqrt{1- \modulus{a}^2 }\right)}^m \int_{ \Space{D}{} }
\frac{f(w)}{{(a\bar{w}-1)}^m} \frac{dw}{{(1- \modulus{w}^2 )}^{2-m}},
\end{eqnarray*}
which for $m=2$ is the classical Bergman formula up to factor 
${\left(\sqrt{1- \modulus{a}^2 }\right)}^m$. Note that calculations in 
standard approaches are ``rather lengthy and must be done in
stages''~\cite[\S~1.4]{Krantz82}. $\eoe$
\end{examplea}
\begin{exampleb}
Now we consider the same group $G=\SL$ but pick up another subgroup 
$H=A$. Let $e_{12}:=e_1e_2$. It follows from~\eqref{eq:p-prop} that 
the mapping from the subgroup $A\sim \Space{R}{}$ to
even numbers\footnote{See Appendix~\ref{pt:bivec-fun} for a definition of 
functions of
even Clifford numbers.} $\chi_\sigma: a \mapsto a^{e_{12}\sigma}= 
{(\exp (e_1e_2 \sigma\ln{a}))}= (a \n{p}_1 + a^{-1} \n{p}_2) ^\sigma$ 
parametrized by $\sigma \in \Space{R}{}$ is a character (in a somewhat
generalized sense). 
It represents an isometric rotation of $\TSpace{T}{}$ by the angle $a$.

Under the present conditions we have from~\eqref{eq:def-s-b} 
and~\eqref{eq:def-r-b}:
\begin{displaymath}
r(g^{-1}*s(\n{u}))=\matr{ \frac{-\n{b}\n{u}+\n{a}}{ 
\modulus{-\n{b}\n{u}+\n{a}} } }{0}{0}{ \frac{-\n{b}\n{u}+\n{a}}{ 
\modulus{-\n{b}\n{u}+\n{a}} } }, \qquad 
g^{-1}=\matr{\n{a}}{\n{b}}{-\n{b}}{\n{a}}.
\end{displaymath}
If we again introduce the exponential coordinates $t$ on $\TSpace{T}{}$ 
coming from the subgroup $A$ (i.e., $\n{u}= e_1 e^{e_1e_2t}\cosh t e_1 -
\sinh t e_2= (x+ \frac{1}{x}) e_1 - (x- \frac{1}{x}) e_2$, $x=e^t$) then
the measure $dt$ on $\TSpace{T}{}$ will satisfy the transformation
condition
\begin{displaymath}
\frac{d(g\cdot t)}{dt}= \frac{1}{(be^{-t}+a)(ce^t+d)}= 
\frac{1}{(-\n{b}\n{u}+\n{a})(\n{u}\n{b}-\n{a})},  
\end{displaymath}
where
\begin{displaymath}
g^{-1}=\matr{a}{b}{c}{d}=\matr{\n{a}}{\n{b}}{-\n{b}}{\n{a}}.
\end{displaymath}

Furthermore we can construct a realization of~\eqref{eq:def-ind} on the 
functions defined on $ \TSpace{T}{} $ by the formula:
\begin{equation} \label{eq:ind-b}
[\pi_\sigma(g) f](\n{v})= 
\frac{(- \n{v} \n{b} + \bar{\n{a}})^{\sigma}}
{(- \n{b} \n{v} + {\n{a}})^{1+\sigma}}
f\left( \frac{\n{a}\n{v}+\n{b}}{-\n{b}\n{v} +\n{a}} \right), \quad 
g^{-1}=\matr{\n{a}}{\n{b}}{-\n{b}}{\n{a}}.
\end{equation}
It is induced by the character $\chi_{\sigma}$ due to formula
$-\n{b} \n{v}+\n{a}= (cx+d)\n{p}_1+(bx^{-1}+a)\n{p}_2$, where $x=e^t$ and 
it is a cousin of the principal series representation (see~\cite[\S~VI.6,
Theorem~8]{Lang85},~\cite[\S~8.2, Theorem~2.2]{MTaylor86} 
and Appendix~\ref{pt:principal}).
The subspaces of vector valued and even number valued functions are invariant 
under~\eqref{eq:ind-b} and the representation is unitary with respect to 
the following inner product (about Clifford valued inner product
see~\cite[\S~3]{Cnops94a}):
\begin{displaymath}
\scalar{f_1}{f_2}_{ \TSpace{T}{} } = \int_{ \TSpace{T}{} }
\bar{f}_2(t) f_1(t)\, dt.
\end{displaymath}
We will denote by $\FSpace{L}{2}(\TSpace{T}{})$ the space of 
$\Space{R}{1,1}$-even Clifford number valued functions on $\TSpace{T}{}$ equipped with
the above inner product.

We select function $f_0(\n{u})\equiv 1$ neglecting the fact that it does 
not belong to $ \FSpace{L}{2}( \TSpace{T}{} )$. Its transformations
\begin{equation} \label{eq:coher-b}
f_g(\n{v})=[\pi_\sigma(g) f_0](\n{v})= \modulus{1+\n{u}^2}^{1/2}
\frac{(- \n{v} \n{b} + \bar{\n{a}})^{\sigma}}{(- \n{b} \n{v} + 
{\n{a}})^{1+\sigma}}
\end{equation}
and in particular the identity $[\pi_\sigma(g)f_0](\n{v})= 
\bar{\n{a}}^\sigma
{\n{a}}^{-1-\sigma} f_0(\n{v})=\n{a}^{-1-2\sigma} f_0(\n{v})$ for 
$g^{-1}=\matr{\n{a}}{0}{0}{\n{a}}$ demonstrates that it is a vacuum vector.
Thus we define the reduced wavelet transform accordingly 
to~\eqref{eq:def-s-b} and~\eqref{eq:def-w} by the formula:
\begin{eqnarray}
[\oper{W}_\sigma f] (\n{u})                        
& = & \int_{ \TSpace{T}{} } \modulus{1+\n{u}^2}^{1/2} \overline{ \left(
\frac{(-e_1e^{e_{12} t} \n{u} + \n{1})^\sigma}
{(-\n{u}e_1 e^{e_{12} t} + \n{1})^{1+\sigma}}
\right) } f(t) \,dt \nonumber \\
& = &  \modulus{1+\n{u}^2}^{1/2} \int_{ \TSpace{T}{} }  
\frac{(-\n{u} e_1 e^{e_{12} t} + \n{1})^\sigma}
{(-e^{-e_{12} t} e_1 \n{u} + \n{1})^{1+\sigma}}
f(t) \,dt \label{eq:sio1} \\
& = &  \modulus{1+\n{u}^2}^{1/2} \int_{ \TSpace{T}{} }  
\frac{(-\n{u} e_1 e^{e_{12} t} +\n{1}  )^\sigma }
{e^{-e_{12} t(1+\sigma)}(- e_1 \n{u} + e^{e_{12} t})^{1+\sigma}}
f(t) \,dt \nonumber \\
& = &  \modulus{1+\n{u}^2}^{1/2} \int_{ \TSpace{T}{} }  
\frac {(-\n{u} e_1 e^{e_{12} t} +\n{1} )^\sigma} 
{(- e_1 \n{u} + e^{e_{12} t})^{1+\sigma}} e^{e_{12} t(1+\sigma)}
f(t) \,dt \nonumber \\
& = &  \modulus{1+\n{u}^2}^{1/2} e_{12}\int_{ \TSpace{T}{} }  
\frac{(-\n{u} e_1 e^{e_{12} t} +\n{1} )^\sigma }
{(- e_1 \n{u} + e^{e_{12} t})^{1+\sigma}}e^{e_{12} t\sigma}
( e_{12} e^{e_{12} t} \,dt )\,
f(t)  \nonumber \\
& = &  \modulus{1+\n{u}^2}^{1/2} e_{12}\int_{ \TSpace{T}{} }  
\frac{(-\n{u} e_1 \n{z}  + \n{1} )^\sigma \n{z}^{\sigma} }
{(- e_1 \n{u} + \n{z})^{1+\sigma}} 
 \,d\n{z} \, f(\n{z})  \label{eq:cauchy-b}
\end{eqnarray}
where $\n{z}=e^{e_{12}t}$ and $d\n{z}=e_{12}e^{e_{12}t}\,dt$ are the new 
monogenic variable and its differential respectively. The 
integral~\eqref{eq:cauchy-b} is a singular one, its four singular 
points are intersections of the light cone with the origin in $\n{u}$ with 
the unit circle $\TSpace{T}{}$. See Appendix~\ref{pt:sio} about the 
meaning of this singular integral operator.

The explicit similarity between~\eqref{eq:cauchy} and~\eqref{eq:cauchy-b}
allows us consider transformation $\oper{W}_\sigma$~\eqref{eq:cauchy-b} as
the analog of the Cauchy integral formula and the linear space $
\FSpace{H}{\sigma}( \TSpace{T}{} ) $~\eqref{eq:hardy-b} generated by the
coherent states $f_\n{u}(\n{z})$~\eqref{eq:coher-b} as the correspondence
of the Hardy space.  Due to ``indiscrete'' (i.e. they are not square 
integrable) nature of principal series representations there are no
counterparts for the Bergman projection and Bergman space.  $\eoe$ 
\end{exampleb}

\subsection{The Dirac (Cauchy-Riemann) and Laplace Operators}
Consideration of Lie groups is hardly possible without consideration of 
their Lie algebras, which are naturally represented by left and right
invariant vectors fields on groups. On a homogeneous space $\Omega=G/H$ we
have also defined a left action of $G$ and can be interested in left
invariant vector fields (first order differential operators). Due to the 
irreducibility of $ \FSpace{F}{2}( \Omega )$ under left action of $G$ every
such vector field $D$ restricted to $ \FSpace{F}{2}( \Omega )$ is a scalar
multiplier of identity $D|_{\FSpace{F}{2}( \Omega )}=cI$. We are
in particular interested in the case $c=0$.
\begin{defn} \cite{AtiyahSchmid80,KnappWallach76}
A $G$-invariant first order differential operator 
\begin{displaymath}
D_\tau: \FSpace{C}{\infty}(\Omega, \mathcal{S} \otimes V_\tau ) 
\rightarrow
\FSpace{C}{\infty}(\Omega, \mathcal{S} \otimes V_\tau )
\end{displaymath}
such that $\oper{W}(\FSpace{F}{2}(X))\subset \object{ker} D_\tau$ is called 
\emph{(Cauchy-Riemann-)Dirac operator} on $\Omega=G/H$
associated with an irreducible representation $ \tau $ of $H$ in a space 
$V_\tau$ and a spinor bundle $\mathcal{S}$.  
\end{defn}
The Dirac operator is explicitly defined by the 
formula~\cite[(3.1)]{KnappWallach76}:
\begin{equation} \label{eq:dirac-def}
D_\tau= \sum_{j=1}^n \rho(Y_j) \otimes c(Y_j) \otimes 1,
\end{equation}
where $Y_j$ is an orthonormal basis of
$\algebra{p}=\algebra{h}^\perp$---the orthogonal completion of the Lie
algebra $\algebra{h}$ of the subgroup $H$ in the Lie algebra $\algebra{g}$
of $G$; $\rho(Y_j)$ is the infinitesimal generator of the right action of
$G$ on $\Omega$; $c(Y_j)$ is Clifford multiplication by $Y_i \in
\algebra{p}$ on the Clifford module $\mathcal{S}$.  We also define 
an invariant Laplacian by the formula
\begin{equation} \label{eq:lap-def}
\Delta_\tau= \sum_{j=1}^n \rho(Y_j)^2 \otimes \epsilon_j \otimes 1,
\end{equation}
where $\epsilon_j = c(Y_j)^2$ is $+1$ or $-1$.
\begin{prop}
Let all commutators of vectors of $ \algebra{h}^\perp $ belong to $
\algebra{h}$, i.e.
$[\algebra{h}^\perp,\algebra{h}^\perp]\subset\algebra{h}$.  Let also $f_0$
be an eigenfunction for all vectors of $ \algebra{h} $ with eigenvalue $0$
and let also $\oper{W}f_0$ be a null solution to the Dirac operator $D$.
Then $\Delta f(x)=0$ for all $f(x)\in \FSpace{F}{2}(\Omega)$.
\end{prop}
\begin{proof}
Because $\Delta$ is a linear operator and $ \FSpace{F}{2}(\Omega)$ is 
generated by $\pi_0(s(a))\oper{W} f_0$ it is enough to check that $
\Delta \pi_0(s(a))\oper{W} f_0=0 $. Because $\Delta$  and $\pi_0$ commute 
it is enough to check that $\Delta \oper{W} f_0=0$. Now we observe that
\begin{displaymath}
\Delta = D^2 - \sum_{i,j} \rho([Y_i,Y_j]) \otimes c(Y_i)c(Y_j) \otimes 1.
\end{displaymath}
Thus the desired assertion is follows from two identities: $D 
\oper{W}f_0=0$ and
$\rho([Y_i,Y_j])\oper{W}f_0=0$, $[Y_i,Y_j]\in H$.
\end{proof}
\begin{examplea}
Let $G=\SL$ and $H$ be its one-dimensional compact subgroup generated by 
an element $Z \in \algebra{sl}(2, \Space{R}{})$ (see 
Appendix~\ref{pt:lie-alg}
for a description of $\algebra{sl}(2, \Space{R}{})$). Then 
$\algebra{h}^\perp$ is spanned by two vectors $Y_1=A$ and $Y_2=B$. In such 
a situation we can use $ \Space{C}{} $ instead of the Clifford algebra $ 
\Cliff[0]{2} $. Then formula~\eqref{eq:dirac-def} takes a simple 
form $D=r(A+iB)$. Infinitesimal action of this operator in the upper-half 
plane follows from calculation in~\cite[VI.5(8), IX.5(3)]{Lang85}, it is 
$[D_{ \Space{H}{} } f] (z)= -2i y \frac{ \partial f(z)}{ \partial \bar{z} }
$, $z=x+iy$. Making the Caley transform we can find its action in the
unit disk $D_{ \Space{D}{} } $: again the Cauchy-Riemann operator $ \frac{ 
\partial }{ \partial \bar{z} } $ is its principal component. 
We calculate  $D_{ \Space{H}{} }$ explicitly now to stress the 
similarity with $\Space{R}{1,1}$ case. 

For the upper half plane $\Space{H}{}$ we have following formulas:
\begin{eqnarray*}
s&:&\Space{H}{} \rightarrow \SL: z=x+iy \mapsto 
g=\matr{y^{1/2}}{xy^{-1/2}}{0}{y^{-1/2}}; \nonumber \\
s^{-1}&:& \SL \rightarrow \Space{H}{}: \matr{a}{b}{c}{d} \mapsto z= 
\frac{ai+b}{ci+d};  \nonumber \\
\rho(g)&:& \Space{H}{} \rightarrow\Space{H}{} :z \mapsto s^{-1}( s(z) * g) 
\\
&& \qquad \qquad \qquad  =s^{-1}\matr{ay^{-1/2}+cxy^{-1/2}}{ 
by^{1/2}+dxy^{-1/2}}{cy^{-1/2}}{dy^{-1/2}}\\
&& \qquad \qquad \qquad =\frac{(yb+xd)+i(ay+cx)}{ci+d} 
\nonumber
\end{eqnarray*}
Thus the right action of $\SL$ on $\Space{H}{}$ is given by the formula
\begin{displaymath}
\rho(g)z=\frac{(yb+xd)+i(ay+cx)}{ci+d}= x+y \frac{bd+ac}{c^2+d^2}+iy 
\frac{1}{c^2+d^2}.
\end{displaymath}
For $A$ and $B$ in $\algebra{sl}(2, \Space{R}{} )$ we have:
\begin{displaymath}
\rho(e^{At}) z = x+iy e^{2t}, \qquad \rho(e^{Bt}) z = x 
+y\frac{e^{2t}-e^{-2t}}{e^{2t}+e^{-2t}}+iy\frac{4}{e^{2t}+e^{-
2t}}.
\end{displaymath}
Thus
\begin{eqnarray*}
[\rho(A) f](z) & = & \frac{ \partial f (\rho(e^{At}) z )}{ \partial t }
|_{t=0} = 2y \partial_2 f(z), \\ {}
[\rho(B) f](z) & = & \frac{ \partial f (\rho(e^{Bt}) z )}{ \partial t } 
|_{t=0} = 2y \partial_1 f(z),
\end{eqnarray*}
where $\partial_1$ and $\partial_2$ are derivatives of $f(z)$ with respect 
to real and imaginary party of $z$ respectively. Thus we get 
\begin{displaymath}
D_{ \Space{H}{} }= i\rho(A) + \rho( B) = 2y i\partial_2 + 2y \partial_1=
2y \frac{ \partial }{ \partial \bar{z} }
\end{displaymath}
as was expected. $\eoe$
\end{examplea}
\begin{exampleb}
In $\Space{R}{1,1}$ the element $B\in \algebra{sl}$ generates the subgroup 
$H$ and its orthogonal completion is spanned by $B$ and $Z$. Thus the
associated Dirac operator has the form $D=e_1 \rho(B) + e_2 \rho(Z)$. 
We need infinitesimal generators of the right action $\rho$ on the ``left'' 
half plane $\TSpace{H}{}$. Again we have a set of formulas similar to the 
classic case:
\begin{eqnarray*}
s&:&\TSpace{H}{} \rightarrow \SL: \n{z}=e_1y+e_2x \mapsto 
g=\matr{y^{1/2}}{xy^{-1/2}}{0}{y^{-1/2}}; \nonumber \\
s^{-1}&:& \SL \rightarrow \TSpace{H}{}: \matr{a}{b}{c}{d} \mapsto \n{z}= 
\frac{ae_1+be_2}{ce_2e_1+d};  \nonumber \\
\rho(g)&:& \TSpace{H}{} \rightarrow \TSpace{H}{} :\n{z} \mapsto s^{-1}( 
s(\n{z}) * g) \\
&&\qquad \qquad \qquad =s^{-1}\matr{ay^{-1/2}+cxy^{-1/2}}{ 
by^{1/2}+dxy^{-1/2}}{cy^{-1/2}}{dy^{-1/2}}\\
&&\qquad \qquad \qquad = \frac{(yb+xd)e_2+(ay+cx)e_1}{ce_2e_1+d}
\nonumber
\end{eqnarray*}
Thus the right action of $\SL$ on $\Space{H}{}$ is given by the formula
\begin{displaymath}
\rho(g)\n{z}=\frac{(yb+xd)e_2+(ay+cx)e_1}{ce_2e_1+d}= 
e_1 y \frac{-1}{c^2-d^2} + e_2 x+e_2y \frac{ac-bd}{c^2-d^2}.
\end{displaymath}
For $A$ and $Z$ in $\algebra{sl}(2, \Space{R}{} )$ we have:
\begin{eqnarray*}
\rho(e^{At}) \n{z}& = & e_1 y e^{2t}+ e_2 x, \\
\rho(e^{Zt}) \n{z}& = &e_1y\frac{-1}{\sin^2 t - \cos^2 t}
+e_2 y\frac{-2\sin t \cos t}{\sin^2 t -\cos^2 t}+ e_2 x\\
& =& e_1y\frac{1}{ \cos 2 t}+ e_2 y \tan 2 t + e_2 x .
\end{eqnarray*}
Thus
\begin{eqnarray*}
[\rho(A) f](\n{z}) & = & \frac{ \partial f (\rho(e^{At}) \n{z} )}{ \partial 
t }
|_{t=0} = 2y \partial_2 f(\n{z}), \\ {}
[\rho(Z) f](\n{z}) & = & \frac{ \partial f (\rho(e^{Zt}) \n{z} )}{ \partial 
t }
|_{t=0} = 2y \partial_1 f(\n{z}),
\end{eqnarray*}
where $\partial_1$ and $\partial_2$ are derivatives of $f(\n{z})$ with 
respect of $e_1$ and $e_2$ components of $\n{z}$ respectively. Thus we get 
\begin{displaymath}
D_{ \TSpace{H}{} }= e_1\rho(Z) + e_2\rho( A) = 2y (e_1\partial_1 + e_2
\partial_2).
\end{displaymath}
In this case the Dirac operator is not elliptic and as a consequence we
have in particular a singular Cauchy integral formula~\eqref{eq:cauchy-b}.
Another manifestation of the same property is that primitives in the
``Taylor expansion'' do not belong to $\FSpace{F}{2}(\TSpace{T}{})$ itself
(see Example~\ref{ex:taylor-b}).  It is known that solutions of a
hyperbolic system (unlike the elliptic one) essentially depend on the
behavior of the boundary value data.  Thus function theory in
$\Space{R}{1,1}$ is not defined only by the hyperbolic Dirac equation alone
but also by an appropriate boundary condition.  
$\eoe$
\end{exampleb}

\subsection{The Taylor expansion}
For any decomposition $f_a(x)=\sum_\alpha  \psi_\alpha(x) V_\alpha(a)$ 
of the coherent states $f_a(x)$ by means of functions $V_\alpha(a)$
(where the sum can become eventually an integral) we have the
\emph{Taylor expansion} 
\begin{eqnarray} 
\widehat{f}(a) & = & \int_X f(x) \bar{f}_a(x)\, dx= \int_X f(x) \sum_\alpha 
\bar{\psi}_\alpha(x)\bar{V}_\alpha(a)\, dx  \nonumber \\
 & = &  \sum_\alpha 
\int_X f(x)\bar{\psi}_\alpha(x)\, dx \bar{V}_\alpha(a) \nonumber \\
 & = & \sum_{\alpha}^{\infty} \bar{V}_\alpha(a) f_\alpha,\label{eq:taylor}
\end{eqnarray}
where $f_\alpha=\int_X f(x)\bar{\psi}_\alpha(x)\, dx$.
However to be useful within the presented scheme such a decomposition 
should be connected with the structures of $G$, $H$, and the representation 
$\pi_0$. We will use a decomposition of $f_a(x)$ by the eigenfunctions of 
the operators $\pi_0(h)$, $h\in \algebra{h}$.
\begin{defn}
 Let $\FSpace{F}{2}=\int_{A} \FSpace{H}{\alpha}\,d\alpha$ be a spectral 
decomposition with respect to the operators $\pi_0(h)$, $h\in \algebra{h}$.
Then the decomposition
\begin{equation} \label{eq:spec-c}
 f_a(x)= \int_{A} V_\alpha(a) f_\alpha(x)\, d\alpha,
\end{equation}
where $f_\alpha(x)\in \FSpace{H}{\alpha}$ and $V_\alpha(a): 
\FSpace{H}{\alpha} \rightarrow \FSpace{H}{\alpha}$ is called the Taylor 
decomposition of the Cauchy kernel $f_a(x)$.
\end{defn}
Note that the Dirac operator $D$ is defined in the terms of left invariant 
shifts and therefor commutes with all $\pi_0(h)$. Thus it also has a 
spectral decomposition over spectral subspaces of $\pi_0(h)$:
\begin{equation} \label{eq:spec-d}
 D= \int_{A} D_\delta \, d\delta.
\end{equation}
We have obvious property
\begin{prop} \label{pr:cauchy-dirac}
If spectral measures $d\alpha$ and $d\delta$ 
from~\eqref{eq:spec-c} and~\eqref{eq:spec-d} have disjoint supports then 
the image of the Cauchy integral belongs to the kernel of the Dirac 
operator.
\end{prop}
For discrete series representation functions $f_\alpha(x)$ can be 
found in $\FSpace{F}{2}$ (as in Example~\ref{ex:taylor-a}), for the 
principal series representation this is not the case. To overcome confusion
one can think about the Fourier transform on the real line. It can be 
regarded as a continuous decomposition of a function $f(x)\in 
\FSpace{L}{2}(\Space{R}{})$ over a set of harmonics $e^{i\xi x}$ neither of
those belongs to $\FSpace{L}{2}(\Space{R}{})$. This has a lot of common 
with the Example~\ref{ex:taylor-b}.
\begin{examplea} \label{ex:taylor-a}
Let $G=\SL$ and $H=K$ be its maximal compact subgroup and $\pi_1$ be 
described by~\eqref{eq:g-transform}.  $H$ acts on $\Space{T}{}$ by
rotations.  It is one dimensional and eigenfunctions of its generator $Z$
are parametrized by integers (due to compactness of $K$).  Moreover, on the
irreducible Hardy space these are positive integers $n=1,2,3\ldots$ and
corresponding eigenfunctions are $f_n(\phi)=e^{i(n-1)\phi}$. Negative 
integers span the space of anti-holomorphic function and the splitting 
reflects the existence of analytic structure given by the Cauchy-Riemann 
equation. The decomposition of coherent states $f_a(\phi)$ by means of this
functions is well known:
\begin{displaymath}
f_a(\phi)= \frac{ \sqrt[]{1- \modulus{a}^2 }}{ \bar{a}e^{i\phi}-1}= 
\sum_{n=1}^\infty \sqrt[]{1- \modulus{a}^2 }\bar{a}^{n-1} e^{i(n-1)\phi}=
\sum_{n=1}^\infty V_n(a)f_n(\phi),
\end{displaymath}
where $V_n(a)=\sqrt[]{1- \modulus{a}^2 }\bar{a}^{n-1} $. This is the 
classical Taylor expansion up to multipliers coming from the invariant 
measure. $\eoe$
\end{examplea}
\begin{exampleb} \label{ex:taylor-b}
Let $G=\SL$, $H=A$, and  $\pi_\sigma$ be described by~\eqref{eq:ind-b}.
Subgroup $H$ acts on $\TSpace{T}{}$ by hyperbolic rotations: 
\begin{displaymath}
\tau: \n{z}=e_1 e^{e_{12}t} \rightarrow e^{2e_{12}\tau}\n{z}=e_1
e^{e_{12}(2\tau+t)}, \qquad t, \tau \in \TSpace{T}{}.
\end{displaymath}
Then for every $p\in \Space{R}{}$ the function 
$f_p(\n{z})=(\n{z})^p=e^{e_{12}pt}$ where $\n{z}=e^{e_{12}t}$ is an
eigenfunction in the representation~\eqref{eq:ind-b} for a generator $a$ 
of the subgroup $H=A$ with the eigenvalue $2(p-\sigma)-1$.  Again, due to
the analytical structure reflected in the Dirac operator, we only need
negative values of $p$ to decompose the Cauchy integral kernel.
\begin{prop} \label{pr:taylor-b}
For $\sigma=0$ the Cauchy integral kernel~\eqref{eq:cauchy-b} has the 
following decomposition:
\begin{equation}
\frac{1}{-e_1\n{u}+\n{z}}=
\int_0^\infty \frac{(e_1\n{u})^{[p]}-1}{e_1\n{u}-1}  
\cdot t\n{z}^{-p}\,dp,
\end{equation}
where $\n{u}=u_1 e_1 +u_2 e_2$, $\n{z}=e^{e_{12}t}$,  and $[p]$
is the integer part of $p$ (i.e. $k=[p]\in \Space{Z}{} $, $k\leq p < k+1$).
\end{prop}
\begin{proof}
Let 
\begin{displaymath}
f(t)=\int_0^\infty F(p) e^{-tp}\, dp
\end{displaymath}
be the Laplace transform. We use the formula~\cite[Laplace Transform Table,
p.~479, (66)]{CRCMTables} 
\begin{equation} \label{eq:lap-table}
\frac{1}{t(e^{kt}-a)}= \int_0^\infty \frac{a^{[p/k]}-1}{a-1} e^{-tp}\, dp
\end{equation}
with the particular value of the parameter $k=1$.
Then using $\n{p}_{1,2}$ defined in~\eqref{eq:p-def} we have
\begin{eqnarray}
\lefteqn{
\int_0^\infty \frac{(e_1\n{u})^{[p]}-1}{e_1\n{u}-1}  
\cdot t\n{z}^{-p}\,dp =}\  && \nonumber \\
& = & t\int_0^\infty \left(\frac{(-u_1-u_2)^{[p]}-1}{(-u_1-u_2)-1} \n{p}_2 
+  \frac{(-u_1+u_2)^{[p]}-1}{(-u_1+u_2)-1} \n{p}_1 
\right)  (e^{tp}\n{p}_2+ e^{-tp}\n{p}_1)\,dp \nonumber \\
& = & t\int_0^\infty  \frac{(-u_1-u_2)^{[p]}-1}{(-u_1-u_2)-1}e^{tp}\,dp\, 
\n{p}_2 +  
t\int_0^\infty \frac{(-u_1+u_2)^{[p]}-1}{(-u_1+u_2)-1}e^{-tp}\,dp\, \n{p}_1 
\nonumber \\
& = & \frac{t}{t(e^{-t} +u_1+u_2)} \n{p}_2 + 
\frac{t}{t(e^{t} +u_1-u_2)} \n{p}_1 \label{eq:lap-appl} \\
& = & \frac{1}{(e^{-t} +u_1+u_2) \n{p}_2 + 
(e^{t} +u_1-u_2) \n{p}_1 }\nonumber \\
&=& \frac{1}{-e_1\n{u}+\n{z}} \nonumber,
\end{eqnarray}
where we obtain~\eqref{eq:lap-appl} by an application 
of~\eqref{eq:lap-table}.
\end{proof}
Thereafter for a function $f(\n{z})\in \FSpace{F}{2}(\TSpace{T}{})$ we 
have the following Taylor expansion of its wavelet transform:
\begin{displaymath}
[\oper{W}_0 f](u)= \int_0^\infty \frac{(e_1\n{u})^{[p]}-1}{e_1\n{u}-1} 
f_p \,dp,
\end{displaymath}
where
\begin{displaymath}
f_p = \int_{\TSpace{T}{}} t \n{z}^{-p}\,d\n{z} f(\n{z}).
\end{displaymath}
The last integral is similar to the Mellin 
transform~\cite[\S~III.3]{Lang85}, \cite[Chap.~8, (3.12)]{MTaylor86}, 
which naturally arises in study of the principal series representations 
of $\SL$. 

I was pointed by Dr.~J.~Cnops that for the Cauchy kernel $(-e_1 
\n{u}+\n{z})$ there is still a decomposition of the form $(-e_1 
\n{u}+\n{z})=\sum_{j=0}^\infty (e_1 \n{u})^j \n{z}^{-j-1}$. It this 
connection one may note that representations $\pi_1$~\eqref{eq:g-transform}
and $\pi_\sigma$~\eqref{eq:ind-b} for $\sigma=0$ are unitary equivalent. 
(this is a meeting point between discrete and principal series). Thus 
a function theory in $\Space{R}{1,1}$ with the value $\sigma=0$ could 
carry many properties known from the complex analysis.
$\eoe$
\end{exampleb}

\setcounter{section}{0}
\renewcommand{\thesection}{\Alph{section}}

\section{Appendix} \label{se:appendix}
\begin{punct}[The Lie algebra of $\SL$] \label{pt:lie-alg}
The Lie algebra $ \algebra{sl}(2,\Space{R}) $ of $\SL$ consists of all 
$2\times 2$ real matrices of trace zero. One can introduce a basis
\begin{displaymath}
A= \frac{1}{2} \matr{-1}{0}{0}{1},\quad B= \frac{1}{2} \matr{0}{1}{1}{0}, 
\quad Z=\matr{0}{1}{-1}{0}.
\end{displaymath}
The commutator relations are
\begin{displaymath}
[Z,A]=2B, \qquad [Z,B]=-2A, \qquad [A,B]=- \frac{1}{2} Z.
\end{displaymath}
\end{punct}
\begin{figure}[h]
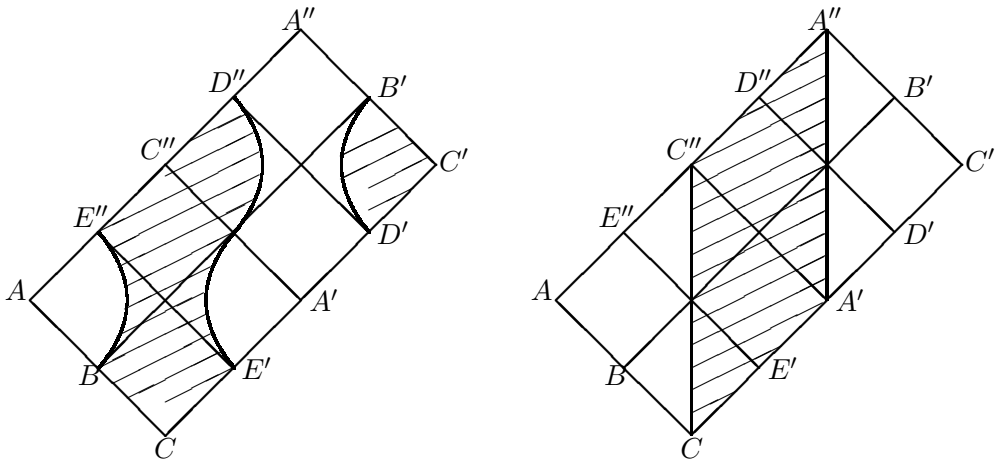

\begin{center}
{\small
\input{r11_a.pic}
\quad
\input{r11_b.pic}
}
\end{center}
\caption{The conformal unit disk, unit circle (on the left) and 
``left'' half plane (on the right) in
$\Space{R}{1,1}$. Points labeled by the same letters should be glued 
together (on each picture separately). Squares $ACA'C''$ and $A'C'A''C''$ 
represent $ \mathbb{R}^{1,1}_- $ and $ \mathbb{R}^{1,1}_+$ respectively. 
Their boundary are the image of the light cone at infinity.} 
\label{pi:disk}
\end{figure}
\begin{punct}[Conformal unit disk, unit circle, and half plane] 
\label{pt:disk}
We are taking two copies $\mathbb{R}^{1,1}_+ $ and $ \mathbb{R}^{1,1}_- $
of $ \Space{R}{1,1} $ glued over their light cones at infinity in such a
way that the construction is invariant under natural action of the M\"obius
transformation.  This aggregate denoted by $ \TSpace{R}{1,1}$ is a two-fold
cover of $ \Space{R}{1,1} $. Topologically $\TSpace{R}{1,1}$ is equivalent
to the Klein bottle. Similar conformally invariant two-fold cover of the 
Minkowski space-time was constructed in~\cite[\S~III.4]{Segal76} in 
connection with the red shift problem in extragalactic astronomy.

We define \emph{(conformal) unit disk} in $ \TSpace{R}{1,1} $ as follows:
\begin{equation}
\TSpace{D}{}=\{\n{u} \such \n{u}^2<-1,\ \n{u}\in \mathbb{R}^{1,1}_+ 
\}
\cup \{\n{u} \such \n{u}^2>-1,\ \n{u}\in \mathbb{R}^{1,1}_- \}.
\end{equation}
It can be shown that $ \TSpace{D}{}$ is conformally invariant and 
has a boundary $\TSpace{T}{}$---the two glued copies of unit circles in $ 
\mathbb{R}^{1,1}_+$ and $\mathbb{R}^{1,1}_-$. 

We call $\TSpace{T}{}$ 
the \emph{(conformal) unit circle} in $\Space{R}{1,1}$. $\TSpace{T}{}$ 
consists of four parts---branches of hyperbola---with subgroup $A\in\SL$ 
acting simply transitively on each of them. Thus we will regard 
$\TSpace{T}{}$ as $ \Space{R}{} \cup \Space{R}{} \cup\Space{R}{} 
\cup\Space{R}{}$ with an exponential mapping $\exp: t \mapsto (+ \mathrm{ 
or } -) e_1^{+ \mathrm{ or } -}$, $e_1^\pm \in \mathbb{R}^{1,1}_\pm$, where
each of four possible sign combinations is realized on a particular copy of
$ \Space{R}{}$. More generally we define a set of concentric circles for 
$-1\leq\lambda<0$:
\begin{equation} \label{eq:circle-l}
\TSpace{T}{\lambda}=\{\n{u} \such \n{u}^2=-\lambda^2,\ \n{u}\in 
\mathbb{R}^{1,1}_+
\}
\cup \{\n{u} \such \n{u}^2=-\lambda^{-2},\ \n{u}\in \mathbb{R}^{1,1}_- \}.
\end{equation}

Figure~\ref{pi:disk} illustrates geometry of the conformal unit disk in 
$\TSpace{R}{1,1}$ as well as the ``left'' half plane conformally equivalent
to it.
\end{punct}
\begin{punct}[Functions of even 
Clifford numbers]
\label{pt:bivec-fun}
Let 
\begin{equation} \label{eq:p-def}
\n{a}= a_1 \n{p}_1 + a_2 \n{p}_2, \qquad \n{p}_1= \frac{1+e_1e_2}{2},\quad 
\n{p}_2= \frac{1-e_1e_2}{2},\quad a_1,a_2\in \Space{R}{}
\end{equation}
be an even Clifford number in $ \Cliff[1]{1} $. It follows from the 
identities
\begin{equation} \label{eq:p-prop}
\n{p}_1 \n{p}_2 = \n{p}_2 \n{p}_1  =0, \qquad \n{p}_1^2=\n{p}_1, \qquad
\n{p}_2^2=\n{p}_2, \qquad \n{p}_1 + \n{p}_2=1
\end{equation}
that $p( \n{a})= p(a_1) \n{p}_1 + p(a_2) \n{p}_2$ for any polynomial
$p(x)$. Let $\FSpace{P}{}$ be a topological space of functions 
$ \Space{R}{} \rightarrow \Space{R}{} $ such that polynomials are dense in 
it. Then for any $ f\in \FSpace{P}{} $ we can define $f(\n{a})$ by the 
formula
\begin{equation} \label{eq:bivec-fun}
f(\n{a})= f(a_1) \n{p}_1 + f(a_2) \n{p}_2.
\end{equation}
This definition gives continuous algebraic homomorphism. 
\end{punct}
\begin{punct}[Principal series representations of $\SL$] 
\label{pt:principal}
We describe a realization of the principal series representations of 
$\SL$. The realization is deduced from the realization by left regular 
representation on the a space of homogeneous function of power $-is-1$ on $ 
\Space{R}{2} $ described in~\cite[\S~8.3]{MTaylor86}. We consider now the 
restriction of homogeneous function not to the unit circle as 
in~\cite[Chap.~8, (3.23)]{MTaylor86} but to the line $x_2=1$ in $ 
\Space{R}{2} $. Then an equivalent unitary representation of $\SL$ acts on 
the Hilbert space $\FSpace{L}{2}( \Space{R}{} )$ with the standard Lebesgue
measure by the transformations:
\begin{equation} \label{eq:principal}
[\pi_{is} (g) f](x)= \frac{1}{ \modulus{cx+d}^{1+is} } f \left( 
\frac{ax+b}{cx+d} \right), \qquad g^{-1}= \matr{a}{b}{c}{d}.
\end{equation}
\end{punct}
\begin{punct}[Boundedness of the Singular Integral Operator 
$\oper{W}_\sigma$] \label{pt:sio}
The kernel of integral operator $\oper{W}_\sigma$~\eqref{eq:sio1} is
singular in four points, which are the intersection of $\TSpace{T}{}$ and
the light cone with the origin in $\n{u}$. One can easily see
\begin{displaymath}
\modulus{\frac{(-\n{u} e_1 e^{e_{12} t} + \n{1})^\sigma}
{(-e^{-e_{12} t} e_1 \n{u} + \n{1})^{1+\sigma}}} =
\modulus{1+\n{u}^2}^{1/2} \frac{1}{ \modulus{t-t_0} }  + O(\frac{1}{ 
\modulus{t-t_0}^2 }).
\end{displaymath}
where $t_0$ is one of four singular points mentioned before for a fixed
$\n{u}$ and $t$ is a point in its neighborhood.  More over the kernel of
integral operator $\oper{W}_\sigma$ is changing the sign while $t$ crossing
the $t_0$.  Thus we can define $\oper{W}_\sigma$ in the sense of the
principal value as the standard singular integral operator.

Such defined integral operator $\oper{W}_\sigma$ becomes a bounded linear
operator $ \FSpace{L}{2} ( \TSpace{T}{} ) \rightarrow \FSpace{L}{2} (
\TSpace{T}{\lambda}) $, where $\TSpace{T}{\lambda}$ is the
circle~\eqref{eq:circle-l} in $\TSpace{R}{1,1}$ with center in the origin
and the ``radius'' $\lambda$.  Moreover the norm of the operator
$\lambda^{-2}\oper{W}_\sigma$ is uniformly bounded for all $\lambda$ and
thus we can consider it as bounded operator
\begin{displaymath}
\FSpace{L}{2} ( \TSpace{T}{} ) \rightarrow \FSpace{H}{\sigma} ( 
\TSpace{D}{} ) ,
\end{displaymath}
where
\begin{equation}
\FSpace{H}{\sigma} ( \TSpace{D}{} ) =\{ f(\n{u}) \such D_{\TSpace{D}{}} f(\n{u} 
)=0,\  \n{u}\in \TSpace{D}{},\
\modulus{\lambda}^{-2}\int_{\TSpace{T}{\lambda}} \modulus{f(\n{u})}^2\, du 
< \infty,\ \forall \lambda<0\}. \label{eq:hardy-b}
\end{equation}
is an analog of the classic Hardy space. Note that 
$\modulus{\lambda}^{-2}d\n{u}$ is exactly the invariant
measure~\eqref{eq:def-m-b} on $\TSpace{D}{}$. 

One can note the similarity of arising divergency and singularities with
the ones arising in quantum field theory. The similarity generated by the 
same mathematical object in basement: a pseudoeuclidean space with 
an indefinite metric.
\end{punct}
\begin{punct}[Open problems] \label{pt:o-problems}
This paper raises more questions than gives answers. Nevertheless it is 
useful to state some open problems explicitly.
\begin{enumerate}
\item Demonstrate that Cauchy formula~\eqref{eq:cauchy-b} is an isometry 
between $\FSpace{F}{2}(\TSpace{T}{})$ and $\FSpace{H}{\sigma}(\TSpace{D}{})$
with suitable norms chosen. This almost follows (up to some constant 
factor) from its property to intertwine two irreducible representations of 
$\SL$.
\item Formula~\eqref{eq:sio1} contains Szeg\"o type kernel, which is domain 
dependent. Integral formula~\eqref{eq:cauchy-b} formulated in terms of 
analytic kernel. Demonstrate using Stocks theorem that~\eqref{eq:cauchy-b} is 
true for other suitable chosen domains.
\item The image of Szeg\"o (or Cauchy) type formulas belong to the kernel 
of Dirac type operator only if they connected by additional condition 
(see Proposition~\ref{pr:cauchy-dirac}). Descriptive
condition for the discrete series can be found 
in~\cite[Theorem~6.1]{KnappWallach76}. Formulate a similar condition for 
principal series representations.
\end{enumerate}
\end{punct}

\renewcommand{\thesection}{}
\section{Acknowledgments}
The paper was written while the author stayed at the Department of
Mathematical Analysis, University of Gent  whose hospitality and support I
gratefully acknowledge. 
The stay was provided by the grant 3GP03196 of the FWO-Vlaanderen (Fund of
Scientific Research-Flanders), Scientific Research Network ``Fundamental
Methods and Technique in Mathematics'' and INTAS grant 93--0322--Ext 
subsequently.

I am grateful to R.~Delanghe, F.~Sommen, and I.~Spitkovski for 
many useful comments and suggestions.  I am in debt to J.~Cnops for his
careful reading of the manuscript and numerous helpful remarks.  This paper 
grew from my talk on ISAAC'97 Congress and I am grateful to R.P.~Gilbert,
J.~Ryan, and D.~Struppa for their support, which allows me attend this
impressing event.
\small
\bibliographystyle{plain}
\bibliography{MRABBREV,analyse,aphysics,arare}
\end{document}